\documentclass[12pt]{article}
\usepackage{latexsym}
\usepackage{amsfonts}
\usepackage{amssymb}
\usepackage{amscd}
\usepackage{amsbsy}
\usepackage{amsmath}
\usepackage[dvips]{graphicx}

\newcommand{\bimn}[7]{\bibitem{#1}#2,\hspace{.25em}{\em #3}, {
#4}\hspace{0.25em}{\bf#5}\hspace{0.25em}(#6)\hspace{0.25em}{#7}.}

\begin{document}
\setlength{\unitlength}{3947sp}
\setlength{\unitlength}{1mm}

 \font\absfont=cmr7
 \font\goth=eufm10
 \def\bfl{{\bf l}}
 \def\bfn{{\bf n}}
 \def\id{{\rm Id}}
 \def\slie{{\goth sl}}
 \def\sltwo{{\slie(2)}}
 \def\glie{{\goth g}}
 \def\ints{{\bf Z}}
 \def\nats{{\bf N}}
 \def\rats{{\bf Q}}
 \def\cx{{\bf C}}
 \def\Vd{{V^*}}
 \def\th{\hbox{$^{\hbox{\absfont th}}$}}
 \def\takes{{\colon}}
 \def\inv{{^{-1}}}
 \def\lra{{\longrightarrow}}
 \def\lmt{{\longmapsto}}
 \def\cgF{{\mathcal F}}
 \def\tr{{\rm tr}}
 \def\binomials#1#2{{\left\{\begin{array}{cc}\kern-4pt#1\kern-4pt\\
  \kern-4pt#2\kern-4pt\\
                             \end{array}\right\}}}
 \def\binomial#1#2{{\left[\begin{array}{cc}\kern-4pt#1\kern-4pt\\
  \kern-4pt#2\kern-4pt\\
                             \end{array}\right]}}
 \def\bin#1#2{{\left(\begin{array}{cc}\kern-4pt#1\kern-4pt\\
  \kern-4pt#2\kern-4pt\\
                             \end{array}\right)}}
 \def\max#1{{{\rm max}\,#1}}
 \def\jsymbol#1#2{{\left({#1\over#2}\right)}}
 \def\chg{\check{c}_\glie}
 %
%


\centerline{\large \bf On Habiro's cyclotomic expansions of the
Ohtsuki invariant} \vskip2ex
 \centerline{Ruth Lawrence and Ofer Ron}
 \centerline{\em Einstein Institute of Mathematics, Hebrew University
of Jerusalem} \centerline{\em Givat Ram, 91904 Jerusalem, ISRAEL}
\centerline{{\em E-mail address:} {\tt ruthel@ma.huji.ac.il}}

%

\begin{abstract}{We give a self-contained treatment of
Le and Habiro's approach to the Jones function of a knot and
Habiro's cyclotomic form of the Ohtsuki invariant for manifolds
obtained by surgery around a knot. On the way we reproduce a state
sum formula of Garoufalidis and Le for the colored Jones function of
a knot. As a corollary, we obtain bounds on the growth of
coefficients in the Ohtsuki series for manifolds obtained by surgery
around a knot, which support the slope conjecture of Jacoby and the
first author.}
\end{abstract}

\section{Introduction} \label{s1}

    Suppose that $M$ is a compact oriented $3$--manifold without
boundary.  For any Lie algebra, $\glie$, and integral level, $k$,
there is defined an invariant, $Z_{k+\chg}(M,L)$, of embeddings of
links $L$ in $M$, known as the Witten-Reshetikhin-Turaev invariant
(see \cite{witten}, \cite{reshetikhinturaev}). It is known that for
links in $S^3$, $Z_K(S^3,L)$ is a polynomial in
$q=\exp{2\pi{}i\over{}K}$, namely the generalized Jones polynomial
of the link $L$.  In this paper we consider only the case
$\glie=su(2)$.

    Now assume that $M$ is a rational homology sphere, with
$H=|H_1(M,\ints)|$.  In the normalization for which the invariant
for $S^3$ is $1$, denote the invariant for the pair $(M,\emptyset)$,
as an algebraic function of $q$ at $K\th$ roots of unity, by
$Z_K(M)$. For a rational homology sphere $M$ and odd prime $K$,
\cite{murakami} showed that $Z_K(M)\in\ints[q]$, so that for some
$a_{m,K}(M)\in\ints$, one has
$$Z_K(M)=\sum_{m=0}^\infty{}a_{m,K}(M)(q-1)^m\>.$$  Although the
$a_{m,K}$ are not uniquely determined, it is known from
\cite{ohtsuki} that there exist rational numbers $\lambda_m(M)$ such
that,
$$a_{m,K}(M)\equiv\lambda_m(M)$$ as elements of $\ints_K$ for all
sufficiently large primes $K$.  For integer homology spheres,
$\lambda_0(M)=1$ and $\lambda_1(M)=6\lambda(M)$ where $\lambda(M)$
denotes the Casson--Walker invariant of $M$ in Casson's
normalization.  As a result, one may define a formal power series
$$Z_\infty(M)=\sum_{m=0}^\infty\lambda_m(q-1)^m\>,$$ with rational
coefficients, which is an invariant of rational homology
$3$--spheres, $M$, known as the {\sl Ohtsuki series} of $M$. This is
the asymptotic expansion of the trivial connection contribution to
$Z_K(M)$ (see \cite{rozansky}). Work of Habiro (see
\cite{habiroone}, \cite{habiro}) showed that $Z_\infty(M)$ can
always be written in the form
$$Z_\infty(M)=\sum_{n=0}^\infty{}a_n(q)(q-1)^n\{n\}!$$
where $a_n(q)\in\ints[q,q^{-1}]$. This expression has the advantage
that at a root of unity it reduces to a finite sum (polynomial in
$q^{\pm1}$, while it also can be considered as an element of
$\ints[[q-1]]$.

In this paper we review the technique used by Le \cite{leint} and
Habiro \cite{habiro}, using the $R$-matrix presentation of link and
manifold invariants from \cite{reshetikhinturaevtwo},
\cite{reshetikhinturaev}, to compute the colored Jones function for
a knot and hence the Ohtsuki series for integer homology spheres
obtained by surgery around a knot.  We carry out the quantum group
calculation in general, obtaining the colored Jones polynomial of a
knot as a state sum directly in terms of $q$-numbers as in
\cite{garoufalidisle} (though here it is written for arbitrary knot
diagrams, not just for braid closures). For arbitrary integer
homology spheres the result would be similar, but the knot case is
simpler, since then the Gauss diagram involves only one circle.

    An outline of the present paper is as follows.  In \S2, basic
notation is defined relating to $q$-numbers along with a summary of
facts needed concerning the quantum group $U_q\slie_2$. In \S3, we
follow \cite{reshetikhinturaev} and Le to give a description of the
$\slie_2$ colored Jones polynomial of a knot, first in terms of
$R$-matrices, then using the universal invariant as an algebra
element and finally evaluating to get the state sum formula
(Theorem~1) for the colored Jones function based on a Gauss diagram
description of the knot.  In \S4, following Le and Habiro, the
transformation from the colored Jones function of a knot to the
Ohtsuki series of the 3-manifolds obtained by surgery around the
knot, is discussed and combining with Theorem~1 gives Habiro's
cyclotomic form for the Ohtsuki series. This formulation is used to
give bounds on the growth of coefficients in the Ohtsuki series,
which support the slope conjecture of \cite{jacobylawrence}. Finally
in \S5, the application of the formulae to the trefoil and figure-8
knots are demonstrated.

\section{$q$-numbers and $U_qsl_2$} \label{qgpsec}

\subsection{$q$-numbers} Let $q=v^2=e^\hbar$ be a formal parameter.
Define $q$-numbers, $q$-factorials and $q$-binomial coefficients
according to
$$[n]={v^n-v^{-n}\over{v-v\inv}}\>,\qquad[n]!=\prod_{i=1}^n{[i]}\>,
 \qquad\binomial{n}{m}={[n]!\over[m]![n-m]!}\>.$$
 The dependence on $q$ will be omitted from the notation. We will
 also define asymmetric $q$-numbers, $q$-factorials and $q$-binomial
 coefficients according to
$$\{n\}={q^n-1\over{q-1}}\>,\qquad\{n\}!=\prod_{i=1}^n{\{i\}}\>,
 \qquad\binomials{n}{m}={\{n\}!\over\{m\}!\{n-m\}!}\>.$$
There is an inductive relation as for ordinary binomial coefficients
$$\binomial{n+1}{r+1}=v^{r-n}\binomial{n}{r}+v^{r+1}\binomial{n}{r+1}\>,\quad
 \binomials{n+1}{r+1}=\binomials{n}{r}+q^{r+1}\binomials{n}{r+1}$$
 from which it follows that the asymmetric $q$-numbers, $q$-factorials
and $q$-binomial coefficients all lie in $\ints[q]$. Also, even when
$x$ is not an integer, say an element of an algebra, we can make
sense of quantum binomial coefficients of the form
$$\binomial{x}{m}=\prod_{i=1}^m{[x-i+1]\over[i]}\>,
\qquad\binomials{x}{m}=\prod_{i=1}^m{\{x-i+1\}\over\{i\}}\>,$$ where
$m\in\nats$. The two sets of $q$-numbers are related by
$$\{n\}=v^{n-1}[n]\>,\qquad\{n\}!=v^{{1\over2}n(n-1)}[n]!\>,
\qquad\binomials{n}{m}=v^{m(n-m)}\binomial{n}{m}\>.$$
 We will also denote by $\{n\}_-$, the asymmetric $q$-number
 obtained by replacing $q$ by $q^{-1}$, and by $\{n\}_-!$ the
 corresponding $q$-factorial.  Thus $\{n\}_-=q^{1-n}\{n\}$ and
 $\{n\}_-!=q^{-{1\over2}n(n-1)}\{n\}!$.

\subsection{The quantum group $U_q\slie(2)$}
 The quantum group $A=U_q\slie(2)$ is defined by generators
$H$, $X$ and $Y$ with relations
$$[X,Y]=[H]\>,\qquad[H,X]=2X\>,\qquad[H,Y]=-2Y\>.$$
 Set $K=e^{\hbar{}H/4}=v^{H/2}$. The comultiplication
 $\Delta\takes{}A\lra{}A\otimes{}A$ is given by
$$
 \Delta(H)=H\otimes1+1\otimes{}H\>,\quad
 \Delta(X)=X\otimes{}K+K^{-1}\otimes{}X\>,\quad
 \Delta(Y)=Y\otimes{}K+K^{-1}\otimes{}Y\>,
  $$
with antipode $S\takes{}A\lra{}A$ being an antihomomorphism
acting on the generators by
$$S(H)=-H\>,\qquad
 S(X)=-vX\>,\qquad
 S(Y)=-v\inv{}Y\>.$$

Additionally $A$ is a quasi-triangular Hopf algebra with universal
$R$-matrix, $R\in{}A\otimes{}A$ satisfying the Yang-Baxter equation
$R_{12}R_{13}R_{23}=R_{23}R_{13}R_{12}$. A formula for $R$ (see
\cite{drinfeld}) is
$$R=\kern-3pt\sum_{l=0}^\infty\kern-4pt{(1-q^{-1})^l\over\{l\}!}
  q^{{1\over4}H\otimes{}H}K^lX^l\otimes{}K^{-l}Y^l
  =\kern-5pt\sum_{l,n=0}^\infty\kern-5pt{(1-q^{-1})^l(\hbar/4)^n\over\{l\}!n!}
  H^nK^lX^l\otimes{}H^nK^{-l}Y^l$$
Write this as $R=\sum_i\alpha_i\otimes\beta_i$ (thus the suffix $i$
refers to a pair $(n,l)$ of non-negative integers). As in any
quasi-triangular Hopf algebra, $R\inv=(\id\otimes{}S\inv)R$,
$$R^{-1}=\sum_{l,n=0}^\infty
{(1-q)^lq^{{1\over2}l(l-1)}(-\hbar/4)^n\over\{l\}!n!}
H^nK^{-l}X^l\otimes{}H^nK^lY^l\>.$$
 We will write these two formulae together for $\sigma=\pm1$ as
 $$R^\sigma=\sum_i\alpha_i^{(\sigma)}\otimes\beta_i^{(\sigma)}
 =\sum_{n,l=0}^\infty{(1-q^{-\sigma})^l
 (\sigma\hbar/4)^n\over\{l\}_\sigma!n!}
 H^nK^{\sigma{}l}X^l\otimes{}H^nK^{-\sigma{}l}Y^l$$
The square of the antipode is then given by
 $$S^2(a)=uau\inv\>,\qquad\hbox{for all $a\in{}A$}\>,$$
 where $u=m(S\otimes\id)R_{21}$ and $uS(u)$ is central.
  Furthermore, $A$ is a ribbon Hopf algebra, that is there is a central (ribbon)
element $v$ such that
 $$v^2=uS(u)\>,\qquad{}S(v)=v\>,\qquad{}\epsilon(v)=1\>,
 \Delta(v)=(R_{12}R_{21})\inv(v\otimes{}v)\>.$$
In our case, $K^2=v\inv{}u$. This special (charmed) element
satisfies the property $S^2(a)=K^2aK^{-2}\quad\forall{}a\in{}A$,
while
$$F\equiv\sum_i\alpha_iK^{-2}\beta_i=\sum_i \beta_iK^2\alpha_i
 =\sum_{l=0}^\infty{(q-1)^l\over\{l\}!}q^{l^2/2}
 e^{\hbar{}H^2/4}K^{-2-2l}X^lY^l\>.$$
lies in the center of $A$.

\subsection{Finite dimensional $U_q\slie(2)$-modules}
 For each $\mu\in\nats$, there is a
 $\mu$-dimensional module $\Lambda_\mu$ with basis
 $\{v_i^{(\mu)}\}_{i=0}^{\mu-1}$.  The quantum group action is described by
 the lowest weight vector $v_0^{(\mu)}$ for which (we omit the
 superscripts $\lambda$)
 $$H(v_0)=(1-\mu)v_0\>,\qquad{}Y(v_0)=0\>,\qquad{}X^i(v_0)=v_i\>.$$
 Since $[X^i,Y]=[i][H-i+1]X^{i-1}$, the whole action is given by
 $$H(v_i)=(2i+1-\mu)v_i\>,\qquad{}X(v_i)=v_{i+1}\>,\qquad
 Y(v_i)=[i][\mu-i]v_{i-1}\>.$$
 Since $f\in{}Z(A)$, it acts as multiplication by a scalar in any
irreducible representation.  In particular, in $\Lambda_\mu$, $f$
acts by scaling by
 $$f_\mu=e^{\hbar(\mu-1)^2/4}v^{\mu-1}=v^{\mu^2-1\over2}\>,$$
as can be verified by direct evaluation on $v_0$, where $H$ and $K$
act as $1-\mu$ and $v^{1-\mu\over2}$ respectively.

\section{The Jones function of a knot}\label{geninvsec}

\subsection{Functorial description}

 Suppose $L$ is an oriented framed link in $S^3$.
 The {\sl generalized colored Jones
polynomial} \cite{reshetikhinturaevtwo} of $L$ is defined when each
component $L_i$ is colored by a representation $V_i$ of $A$ and will
be denoted by $J_L(V_1,\ldots,V_c)$. We will assume for simplicity
that $V_i$ are irreducible.  Indeed, according to
\cite{reshetikhinturaevtwo}, there is a functor, $\cgF$, from the
category of (colored) ribbon tangles to the category of vector
spaces, under which links (closed tangles) map to scalars, namely
$J_L$.  A slice of a colored (oriented) ribbon tangle (object in the
category) is an ordered list of colors $V_1,\ldots,V_r$ with
orientations $\epsilon_1,\ldots,\epsilon_r$ and is mapped under the
functor $\cgF$ to
$V_1^{\epsilon_1}\otimes\cdots\otimes{}V_r^{\epsilon_r}$, where
$V^-$ denotes the dual representation $\Vd$.  This functor is
defined by its images on the generators shown in Figure~1. In our
conventions, downward oriented strands are counted positively in
objects, while tangle morphisms are composed considering them from
bottom to top.

\begin{figure}
  \includegraphics[width=350pt]{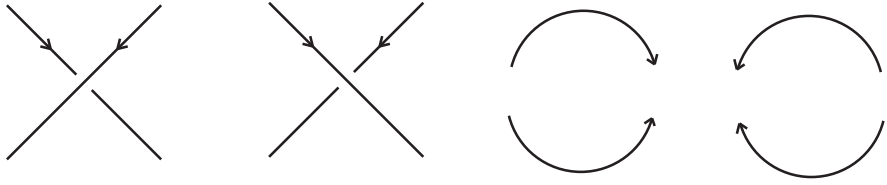}\\
  \caption{Generators of the category of tangles}\label{tanglegen}
\end{figure}

 The two orientations of crossing, in which the strands
are labeled with representations $V$ and $W$ as shown, have as
images the maps $V\otimes{}W\lra{}W\otimes{}V$ given by $P\circ{}R$
and $R\inv\circ{}P$, respectively, where $P$ is the permutation of
the factors.  The four cup and cap sections transform to maps
 \begin{eqnarray*}
 {\Vd}\otimes{}V &\lra& \cx\qquad\qquad{}V\otimes{\Vd}\lra\cx\\
 (x,y) &\lmt& {}x(y)\qquad\qquad(y,x)\lmt{}x(K^2y)
 \end{eqnarray*}
 and
 \begin{eqnarray*}
 \cx &\lra& V\otimes{\Vd}\qquad\qquad{}\cx\lra{\Vd}\otimes{}V\\
 1 &\lmt&
 \sum_ie_i\otimes{}e^i\qquad\quad1\lmt\sum_iK^{-2}e^i\otimes{}e_i
 \end{eqnarray*}
 respectively. Here $\{e^i\}$ is the dual basis for $\Vd$ to the
 basis $\{e_i\}$ for $V$.

Suppose that $L$ has one component (knot) and $T$ is a 1-tangle
presentation of $L$, that is $T\in{\rm Morph}(a,a)$ (where $a$ is
the object consisting of one downward oriented point) such that its
closure with blackboard framing is ambient isotopic to $L$.  Then
 the above
prescription can be used to compute $\cgF(T)$ (when $a$ and $T$ are
colored by $V$) as a map $V\lra{}V$, for any representation $V$ (the
color of the one open strand).  This map commutes with the action of
$A$, so that if $V$ is irreducible it is given by multiplication by
a scalar, namely
 $$J'_L(\mu)={J_L(\Lambda_\mu)\over{}J_U(\Lambda_\mu)}\>,$$
 where $U$ is the unknot with framing zero.  For $U$ we have
$$J_U(\Lambda_\mu)=\tr_{\Lambda_\mu}(K^2)=\tr_{\Lambda_\mu}(K^{-2})
=v^{\mu-1}+v^{\mu-3}+\cdots+v^{1-\mu}=[\mu]\>.$$
 By construction, $J'_L$ is multiplicative under connect sum on $L$
(corresponding to composition of tangles $T$).  In particular, the
algebra element corresponding to a trivial 1-tangle with framing 1
is $f$, so that $J'_L$ changes by a factor $(f_\mu)^t$ under a
framing change of $t$ in $L$.

\subsection{Universal $\sltwo$ formulation of Jones polynomial of a knot}

The above functorial prescription for $J'_L$ when $L$ is a knot
expressed as the closure of a 1-tangle $T$, may be rewritten
algebraically as follows. First assume that $T$ is presented as a
tangle diagram in generic position, in which the crossings are
between downward oriented strands (this is always possible). Recall
that $R^\sigma=\sum_i\alpha^{(\sigma)}_i\otimes\beta^{(\sigma)}_i$
where the sum is over $i\in{}I$ an indexing set (in this case, pairs
of non-negative integers $(n,l)$).  Place an element of $I$ at each
crossing.  In the neighborhood of each crossing decorate the
crossing stands by elements of $A$ according to Figure~2;
$\alpha^\sigma_i$ on the overcrossing arc and $\beta^\sigma_i$ on
the undercrossing arc, where $\sigma$ denotes the sign of the
crssing ($+$ for a slash crossing and $-$ for a backslash crossing).
At local maxima/minima oriented leftwards decorate the strand by
$K^2$, $K^{-2}$ respectively; rightward oriented cups and caps
receive no decoration.

\begin{figure}
  \includegraphics[width=400pt]{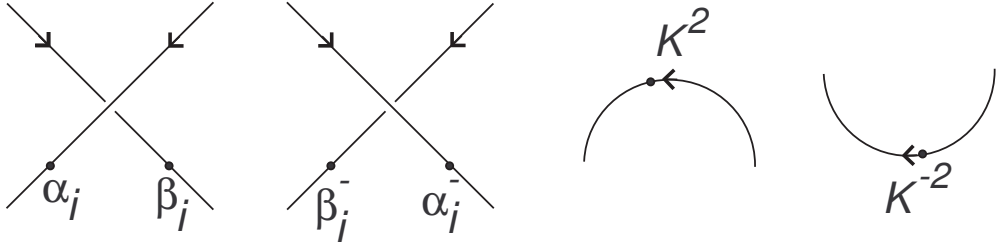}\\
  \caption{Decorations on 1-tangle elements}\label{tangledec}
\end{figure}

 Then ${\cgF}(T)$ is the evaluation in the chosen representation of
the algebra element read off the diagram by tracing the tangle
strand according to its orientation, and writing down the
decorations in the form of a product from left to right, and then
summing over all labels $i$ at crossings. In fact, in the correct
quotient of $A$, the algebra element itself is an invariant (this is
the universal $U_h\slie(2)$ invariant); however we do not need this
here.

In a general tangle diagram, crossings need not always occur only
between downward oriented strands.  However, other crossings can
always be redrawn in terms of such downward oriented crossings, at
the cost of introducing extra caps and cups, by ``internally
rotating" the crossing. For example a crossing of two upward
oriented strands can be obtained by ``internally rotating" a
downward oriented crossing through $\pi$, either clockwise or
anticlockwise and consistency of the results follows from the fact
that $(S^2\otimes{}S^2)R=R$ along with $S^2(a)=K^2aK^{-2}$
$\forall{}a\in{}A$.

\begin{figure}
  \includegraphics[width=400pt]{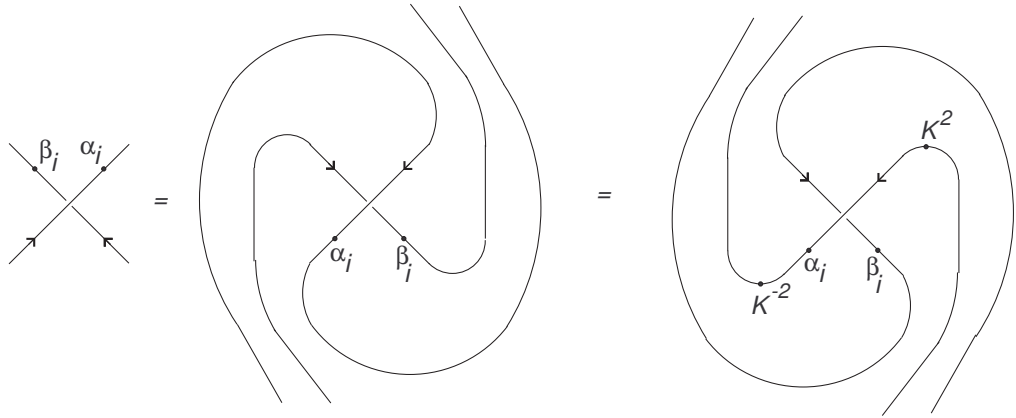}\\
  \caption{Internal rotation of a crossing}\label{crossingrot}
\end{figure}

There are altogether eight crossing types given by the orientations
on the two strands as well the sign of the crossing itself.   The
result is that for an arbitrarily oriented crossing, the
undercrossing arc is decorated with $\beta_i^{(\sigma)}$ while the
overcrossing arc is decorated with $\alpha_i^{(\sigma)}$ {\it
except} in the case of a left pointing crossing in which case it is
decorated by $S^{2\sigma}(\alpha_i^{(\sigma)})$.

\begin{figure}
  \includegraphics[width=400pt]{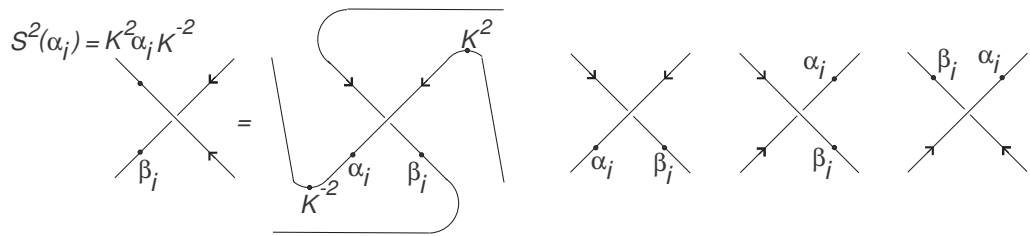}\\
  \caption{Decorations on all orientations of a positive crossing}
  \label{positivecrdec}
\end{figure}

\subsection{Gauss diagram formulation of Jones function of a knot}

We will now proceed to carry out the above prescription at the
quantum group level for an arbitrary 1-tangle, producing a
prescription for $J'_L(\mu)$ as a state sum of scalar quantities.

The combinatorial data involved in specifying a knot $L$, presented
as a 1-tangle closure can be contained in its Gauss diagram.  That
is, consider the knot as parametrised by an (orieted) circle. Each
crossing in $L$ corresponds to exactly two points on the circle, one
for the underpass and the other for the overpass. Encode this by
joining the points on the circle by an oriented chord from overpass
to underpass, marked by a sign to denote the sign of the crossing.
The resulting diagram is known as the {\it Gauss diagram} of $L$; it
is an oriented signed chord diagram on an oriented circle. This is
now sufficient data to encode the knot.

Unfortunately the above description of the algorithm for $J'_L(\mu)$
requires also a tangle presentation with left pointing
crossings/cups/caps marked.  This requires that the circle be marked
with a basepoint as well as blobs $K^{\pm2}$ at left pointing
cups/caps with additionally left pointing crossings marked. (Note
that only left pointing cups/caps/crossings in the long knot are
marked; these may differ from those in the closure). Since
$S^{2\sigma}(a)=K^{2\sigma}aK^{-2\sigma}$, instead of marking a
chord as ``left pointing" (whose corresponding crossing was left
pointing), one may equivalently place $K^{2\sigma}$ and
$K^{-2\sigma}$ blobs on the outer circle just before and after the
outward pointing (corresponding to the overcrossing arc) end of that
chord. The result is an {\it enhanced Gauss diagram} in which
additional marks have been placed on the outer circle, namely the
basepoint and blobs $K^{\pm2}$.  Figure~5 shows an example for the
figure 8 knot.

\begin{figure}
  \includegraphics[width=400pt]{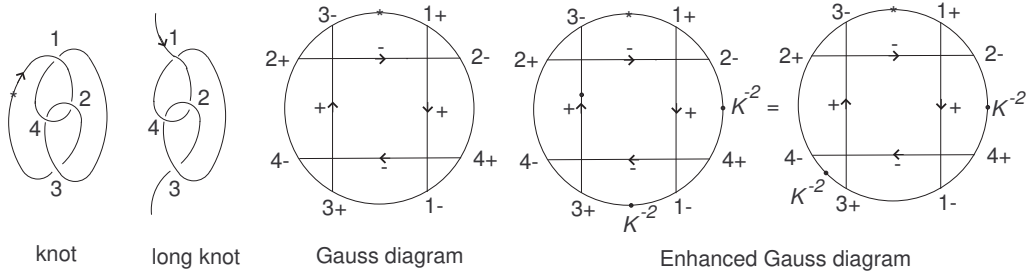}\\
  \caption{Gauss diagrams and enhanced Gauss diagrams}
  \label{gaussdiag}
\end{figure}

Let $D$ be such an enhanced Gauss diagram for a knot $K$ with $c$
crossings. In order to write formulae more easily, we give complete
combinatorial data for $D$. Number the crossings (chords of $D$) by
$j=1,\ldots,c$. Let $\sigma_j$ denote the sign (orientation) of the
crossing, positive for a slash crossing pointing up/down. As we go
around the circle, we encounter $2c$ contributions from crossings
which we will index by $i=1,\ldots,2c$; for each such $i$ there is a
pair $(j(i),\epsilon_i)$ giving the crossing number $j=j(i)$ and the
sign with which it is encountered, $\epsilon_i$, which is positive
for an overcrossing and negative for an undercrossing. Finally let
$i_1,\ldots,i_b$ denote the positions of the blobs $K^{\pm2}$, the
$a\th$ blob appearing between the $i_a\th$ and $(i_a+1)\th$
positions on the circle, with contribution $K^{2\delta_a}$.
Throughout the remainder of this paper, we will always use $j$ as an
index for crossings and $i$, $k$ as indices for vertices on the
Gauss diagram.

The evaluation of $\cgF(T)$ in the representation $\Lambda_\lambda$
may be carried out by computing the action on the lowest weight
vector $v_0^{(\lambda)}$ (equivalently it is the $00$ matrix
element), and then
 \begin{eqnarray}\label{jonesqg}
J'_L(\mu)=\sum_{n_j,l_j=0\atop{}j=1,\ldots,c}^\infty
\left(\prod_{j=1}^c{(1-q^{-\sigma_j})^{l_j}
 (\sigma_j\hbar/4)^{n_j}\over\{l_j\}_{\sigma_j}!.n_j!}\right)
 W^\mu_D(\bfn,\bfl)\>.
 \end{eqnarray}
  Here $W^\mu_D(\bfn,\bfl)$ is the $00$ matrix element in the
evaluation on $\Lambda_\mu$ of the quantum group word
$W_D(\bfn,\bfl)$ obtained by tracing round the outer circle of $D$
starting at the basepoint, with contributions of
$H^{n_j}K^{\sigma_j{}l_j}X^{l_j}$ (outward arrow/overcrossing),
$H^{n_j}K^{-\sigma_j{}l_j}Y^{l_j}$ (inward arrow/undercrossing) and
$K^{\pm2}$ (blobs on the circle), written in order from left to
right. It is convenient to use the notation $X^+=X$ and $X^-=Y$.
Using
$$(X^\epsilon)^lH=(H-2\epsilon{}l)(X^\epsilon)^l\>,\quad{}
(X^\epsilon)^lK=v^{-\epsilon{}l}K(X^\epsilon)^l\>,$$
 we push all powers of $X$ and $Y$ to the end of the word, obtaining
$W_D(\bfn,\bfl)=A(\bfn,\bfl)B(\bfl)$ where
 \begin{eqnarray}\label{bl}
 B(\bfl)=\prod_{i=1}^{2l}(X^{\epsilon_i})^{l_{j(i)}}\>,
 \end{eqnarray}
the product being written from left to right. (Note that $B(\bfl)$
is independent of $n_1,\ldots,n_c$.) The term $A(\bfn,\bfl)$ is the
product of $H$'s and $K's$ remaining after all the powers of $X$ and
$Y$ have been pushed to the right, namely
 $$A(\bfn,\bfl)
 =\bigg[\prod_{i=1}^{2c}\big(H-2\sum_{k<i}\epsilon_kl_k\big)^{n_i}
 \big(\prod_{k<i}v^{-\epsilon_kl_k.\epsilon_i\sigma_il_i}\big)
 K^{\epsilon_i\sigma_il_i}\bigg]\prod_{a=1}^b
 \big(\prod_{k<i_a}v^{-\epsilon_kl_k.2\delta_a}\big) K^{2\delta_a}\>.
 $$
 where by abuse of notation we have written $l_i$ and $n_i$ in place of $l_{j(i)}$
and $n_{j(i)}$ respectively. Collecting terms and noting that
$\sum_i\epsilon_i\sigma_il_i=0$,
$$A(\bfn,\bfl)
=v^{-\sum_i\sum_{k<i}\epsilon_i\epsilon_k\sigma_i{}l_il_k
-2\sum_a\sum_{k<i_a}\delta_a\epsilon_kl_k}K^{2\sum_a\delta_a}
\prod_{i=1}^{2c}\big(H-2\sum_{k<i}\epsilon_k{}l_k\big)^{n_i} \>.$$
 Since both $A(\bfn,\bfl)$ and $B(\bfl)$ preserve the weight,
$W^\mu_D(\bfn,\bfl)=A^\mu(\bfn,\bfl)B^\mu(\bfl)$ where the
superscripts $\mu$ denote the $00$ matrix element in the evaluation
on $\Lambda_\mu$. Thus $A^\mu(\bfn,\bfl)$ is obtained from the above
expression by replacing all occurrences of $H$ and $K$ by $1-\mu$
and $v^{1-\mu\over2}$, respectively.  We may now carry out the sum
over $n_1,\ldots,n_c$ in (\ref{jonesqg}) for $J'_L(\mu)$, since the
sum on $n_j$ may be taken out to give a factor (where
$\lambda=\mu-1$)
$$\bigg[\kern-1pt\sum_{n_j=0}^\infty{(\sigma_j\hbar/4)^{n_j}\over{}n_j!}
\big(H-2\kern-7pt\sum_{k<j+}\kern-7pt\epsilon_k{}l_k\big)^{n_j}
\big(H-2\kern-7pt\sum_{k<j-}\kern-7pt\epsilon_k{}l_k\big)^{n_j}\bigg]_{00}
 \kern-14pt=v^{{\sigma_j\over2}(\lambda+2\sum\limits_{k<j+}\epsilon_k{}l_k)
(\lambda+2\sum\limits_{k<j-}\epsilon_k{}l_k)}$$
 Moving to the zero framed knot introduces a framing correction of
$f_\mu^{-\sum_j\sigma_j}$, that is
$\prod_jv^{-{\sigma_j\over2}(\mu^2-1)}$. The contribution of the
term $A^\mu(\bfn,\bfl)$ to $J'_{L_0}(\mu)$ becomes simply a power of
$v$.  Substituting into (\ref{jonesqg}), we get the zero framed
colored Jones polynomial as
 \begin{eqnarray}\label{jonesq}
J'_{L_0}(\mu)=\sum_{l_j=0\atop{}j=1,\ldots,c}^\infty
\left(\prod_{j=1}^c{(1-q^{-\sigma_j})^{l_j}
\over\{l_j\}_{\sigma_j}!}\right)
 v^{a^\mu(\bfl)}B^\mu(\bfl)\>,
 \end{eqnarray}
where (again $\lambda=\mu-1$)
\begin{eqnarray*}
a^\mu(\bfl)= &-&\sum_i\sum_{k<i}\epsilon_i\epsilon_k\sigma_i{}l_il_k
-2\sum_a\sum_{k<i_a}\delta_a\epsilon_kl_k-\lambda\sum_a\delta_a\\
 &+&\sum_j\sigma_j
 \bigg(-\lambda+\lambda\sum_{k<j-}\epsilon_kl_k
       +\lambda\sum_{k<j+}\epsilon_kl_k
       +2\big(\sum_{k<j-}\epsilon_kl_k\big)
         \big(\sum_{k<j+}\epsilon_kl_k\big)\bigg)
\end{eqnarray*}
 The function $a^\mu(\bfl)$ is a linear function of $\mu$, and
depends quadratically on the parameters $l_j$. Looking at only the
parity of terms (that is, mapping onto $\ints/2\ints[\mu])$,
$$a^\mu(\bfl)\equiv
\sum_i\sum_{k<i}\epsilon_i\epsilon_kl_il_k
+\big(b+c+\sum_i\sum_{k<i}l_k\big)(\mu-1)
$$
Since $\sum_i\epsilon_il_i=0$ (each term $l_j$ appears twice with
opposite signs), the first term is $-\sum_jl_j^2$. Reversing the
order of summation in the last term,
$\sum_i\sum_{k<i}l_k\equiv\sum_i{il_i}\equiv\sum_jp_jl_j$, where
$p_j$ is the number of segments on the Gauss diagram between the two
occurrences of the crossing $j$, that is between $j+$ and $j-$.
Observing that $b+c$ is always even and that $p_j$ is always odd (as
follow from changing the signs of crossings in the knot diagram so
as to obtain an unknot and then applying Reidemeister moves to untie
it and observing that in the process these parities are preserved),
it follows that
$$a^\mu(\bfl)\equiv\mu\sum_jl_j\quad({\rm mod}\,2)\>.$$

The final evaluation required in order to remove all traces of the
quantum group $A$ from our formula for $J'_{L_0}(\mu)$ is that of
$B^\mu(\bfl)$ in (\ref{bl}), namely the $00$ matrix element in the
representation $\Lambda_\mu$ of the product of terms $X^{l_i}$ and
$Y^{l_i}$ described by the Gauss diagram. Using the explicit form
for the basis for $\Lambda_\mu$, one can read the product $B(\bfl)$
from right to left, each term $X$ pushes up the level, while each
$Y$ pushes it down with a factor $[i][\mu-i]$ starting at level $i$.
Reading instead from left to right, $X$ decreases the level while
$Y$ raises the level and introduces the factor $[i+1][\mu-1-i]$ when
starting from level $i$. Thus we obtain
$$B^\mu(\bfl)=\prod_{j=1}^c
\prod_{i=s(j)}^{s(j)+l_j-1}[i+1][\mu-1-i]$$
 where $s(j)=-\sum\limits_{i<j-}l_i\epsilon_i$ is the level at the
start of the term $Y^{l_j}$ contributed by $j-$ in the Gauss diagram
(if any $s(j)$ is negative, the expression vanishes).  Writing
$[i+1][\mu-1-i]=v^{-\sigma\mu}\{i+1\}_\sigma\{\mu-1-i\}_\sigma$, we
obtain
$$B^\mu(\bfl)=\prod_{j=1}^cv^{-\mu{}l_j\sigma_j}
{\{s(j)+l_j\}_{\sigma_j}!\over\{s(j)\}_{\sigma_j}!}
\{\mu-s(j)-1\}_{\sigma_j}\cdots\{\mu-s(j)-l_j\}_{\sigma_j}\>.
$$
Substituting into (\ref{jonesq}) and putting
$c^\mu(\bfl)={1\over2}a^\mu(\bfl)
-{1\over2}\mu\sum_j\sigma_jl_j\in\ints[\mu,l_j]$, we obtain the
following result.\vskip2ex

\noindent{\bf Theorem 1 \cite{garoufalidisle}}{\em \quad The colored
Jones function of a zero framed knot is given by
$$J'_{L_0}(\mu)=\sum_{l_j=0\atop{}j=1,\ldots,c}^\infty
q^{c^\mu(\bfl)} \prod_{j=1}^c \binomials{s(j)+l_j}{l_j}_{\sigma_j}
 (q^{\sigma_j(\mu-s(j)-1)}-1)\cdots
    (q^{\sigma_j(\mu-s(j)-l_j)}-1)\>,
$$
where $s(j)=-\sum_{i<j-}\epsilon_il_{j(i)}$ and $c^\mu(\bfl)$ is a
function dependent on the Gauss diagram, linear in $\mu$ and
quadratic in $l_j$, with integer coefficients, given explicitly by
(\ref{dl}) and (\ref{cl}) below.}\vskip2ex

  Several remarks and standard
facts about the colored Jones function follow immediately from this
theorem.

\begin{itemize}
  \item  In \cite{garoufalidisle}, this theorem is only stated
explicitly for braid closures. A very similar (though inequivalent)
state sum for the colored Jones function was recently obtained in
Theorem 7, \cite{garoufalidisloebl}.
\item Considering $q=1+h$, each term can be thought of as a power
series in $h$ whose coefficients are polynomials in $\mu$.  The
contribution from a particular set of $c$ non-negative integers
$l_j$ is divisible by $h^{\sum_jl_j}$.  That is, the coefficient of
$h^N$ comes only from $l_j$'s with $\sum_jl_j\leq{}N$, a finite
number of contributions. So $J'_{L_0}(\mu)$ is meaningful as a
formal power series in $h$ whose coefficients are polynomials in
$\mu$ (even when $\mu$ is not an integer); this is the Jones
function.
 \item If $\mu\in\nats$, then the contribution from any term
for which $s(j)\geq\mu-1$ will vanish. Thus only a finite number of
$\bfl$'s will contribute, giving a polynomial; this is the colored
Jones polynomial.
  \item For each $l$, the contribution is a polynomial in $q$ and
 $q^\lambda$, so that $J'_{L_0}(\lambda)$ can be written as a power
 series in $h$ whose coefficients are polynomials in $q^\lambda$.
 (See the Melvin-Morton-Rozansky conjecture of \cite{melvinmorton},
 \cite{rozanskymmr} proved in \cite{barnatanmmr}.)
 \item Since
$\binomials{n}{m}\in\ints[q]$, it follows that the contribution from
$l$ is divisible by $h^{\sum_jl_j}\prod_j\{j\}!$.  At a root of
unity of order $K$, only terms with $l_j<K$ $\forall{}j$ will
contribute. Again there are a finite number of contributions.
 \item The coefficient of $\mu$ in
$c^\mu(\bfl)$ simplifies to
 \begin{eqnarray}\label{dl}
 d(\bfl)=\sum_jq_jl_j-{1\over2}\sum_a\delta_a-{1\over2}\sum_j\sigma_j
 \end{eqnarray}
 where $q_j={1\over2}\sum_{i>j+}\sigma_i
 -{1\over2}\sum_{i>j-}\sigma_i-{1\over2}\sigma_j\in\ints$. When the
knot is presented in a blackboard framing which is also zero framed,
then $\sum_j\sigma_j=0$, so that
$$d(\bfl)=\sum_jq_jl_j-{1\over2}\sum_a\delta_a$$ where now $q_j$ has
the meaning of half the sum of the signs of crossings encountered
when tracing around the knot from the overcrossing arc of the $j\th$
crossing to its undercrossing arc; equivalently, it is the linking
number of the two parts of the knot which the $j\th$ crossing, if
split so as to give a link, separates.
 \item In the limit $q\lra1$ with $q^\lambda$ (or $q^{\lambda+1}$) constant, we obtain a
 formal power series in $h=q-1$.  In particular, the Kashaev
 invariant obtained by looking at the colored Jones polynomial in
 color $K$ at the $K\th$ root of unity $q=e^{2\pi{}i\over{}K}$ is
 obtained by setting $\mu=K$, $q^\mu=1$,
 $$K_L(q)=\kern-13pt\sum_{l_j=0\atop{}j=1,\ldots,c}^\infty
\kern-10pt{}q^{c^0(\bfl)}
\prod_{j=1}^c(q^{-\sigma_j}-1)\ldots(q^{-\sigma_j{}l_j}-1)
\binomials{s(j)+l_j}{l_j}_+\binomials{s(j)+l_j}{l_j}_-
$$
where $c^0(\bfl)$ is the constant term in
$c^{\mu}(\bfl)=d(\bfl)\mu+c^0(\bfl)$, namely
 \begin{eqnarray}\label{cl}
 c^0(\bfl)=-d(\bfl)
 -\kern-6pt\sum_{a\atop{k<i_a}}\kern-6pt\delta_a\epsilon_kl_k
 -{1\over2}\sum_j\kern-6pt\sigma_jl_j
 -{1\over2}\sum_{i\atop{k<i}}\sigma_i\epsilon_i\epsilon_kl_il_k
 +\kern-10pt\sum_{j\atop{i<j-,k<j+}}\kern-10pt
        \sigma_j\epsilon_i\epsilon_kl_il_k
 \end{eqnarray}

\item When $\mu=2$, the colored Jones polynomial reduces to the
one-variable Jones polynomial \cite{jones}.  It can be checked that
the only non-zero contributions to the sum in Theorem~1 in this case
come from $l_j\in\{0,1\}$ for which
$-\sum_{k<i}\epsilon_kl_k\in\{0,1\}$ for all $i$.  This sum is
similar to (though different from) the beautifully simple state
model given by Kauffman for the bracket polynomial in
\cite{kauffman}.
 \end{itemize}

Put concisely, the colored Jones function depends on two parameters,
$q$ and $\mu$, and one can make sense of it if at least one of the
following holds: (i) $q$ is a root of unity; (ii) $\mu\in\nats$;
(iii) $q=1+h$ is considered as a formal parameter with expressions
as formal power series in $h$ and $\mu$ is fixed; (iv) $q$ is
considered as a formal parameter with expressions as formal power
series in $h$ and $q^{\mu}$ is fixed.

\section{The WRT and Ohtsuki invariants}

\subsection{The $SU(2)$ and $SO(3)$ 3-manifold invariants}
Suppose that $q$ is a $K\th$ root of unity.  The prescription of
\cite{reshetikhinturaev}  defines for a link $L$ with $c$ components
and framing $f_i=\pm1$ on the $i\th$ component,
 $$\langle{}L\rangle=\sum_{\mu_1,\ldots,\mu_c=1}^{K-1}
 J'_{L_0}(\mu_1,\ldots,\mu_c)
 \prod_{i=1}^c[\mu_i]^2q^{\sum_if_i\Omega(\mu_i)}\>,$$
 where $q^{\Omega(\mu)}=f_\mu$ is the framing normalization,
that is $\Omega(\mu)={\mu^2-1\over4}$.

  The $SU(2)$ invariant of a 3-manifold $M$ obtained from $S^3$ by
surgery around a link $L$ is then
 $$Z_K(M)=G_+^{-\sigma_+}G_-^{-\sigma_-}\langle{}L\rangle\>,$$
where $\sigma_\pm$ are the numbers of positive/negative eigenvalues
of the linking matrix of $L$ and $G_{\pm}$ are the bracket values
for the unknot with framings $\pm1$.  When $L$ is a knot ($c=1$)
with framing $f=\pm1$, so that $M$ is an integer homology 3-sphere,
this reduces to
 $$Z_K(M)={\sum\limits_{\mu=1}^{K-1}q^{f\mu^2/4}[\mu]^2
 J'_{L_0}(\mu)
 \over\sum\limits_{\mu=1}^{K-1}q^{f\mu^2/4}[\mu]^2}\>.$$

The $SO(3)$ invariant is obtained by the same procedure except that
the sums are all restricted to $\mu_i$ even; see \cite{kirbymelvin}.

\subsection{The $SU(2)$ perturbative invariant (Ohtsuki series)}

Following Le and Habiro, we now show how to obtain Habiro's
cyclotomic form for the Ohtsuki series from Theorem~1.  Assume that
$M$ is obtained by $f=\pm1$-surgery around a knot $L$, so that $M$
is an integer homology sphere.

The expression for $J'_{L_0}(\mu)$ in Theorem~1 is a sum over $l$ of
terms which are polynomial in $q^\mu$ and $q$ (and their inverses).
A term $q^{m\mu}$ in $J'_{L_0}(\mu)$ contributes
$${\sum\limits_{\mu=1}^{K-1}q^{f\mu^2/4}[\mu]^2q^{m\mu}
 \over\sum\limits_{\mu=1}^{K-1}q^{f\mu^2/4}[\mu]^2}$$
to $Z_K(M)$. Both summands vanish at $\mu=0$ and are even in $\mu$,
so that the sum can be replaced by a complete period $2K$.
Completing the square,
$\sum{}q^{f\mu^2/4}q^{m\mu}=q^{-fm^2}\sum{}q^{f\mu^2/4}$, thus the
contribution to $Z_K(M)$ of a term $q^{m\mu}$ in $V'_{L_0}(\mu)$ is
$${q^{-f(m+1)^2}+q^{-f(m-1)^2}-2q^{-fm^2}\over{2q^{-f}-2}}
\in\ints[q^{-f}]\>.$$

Let $B$ be the ring of formal power series in $h=q-1$ whose
coefficients are polynomials in $q^{\pm\mu}$ with integer
coefficients. For $f=\pm1$, define a map
$\phi_f\takes{}B\lra\ints[[h]]$ linear over $\ints[[h]]$ by
$$\phi_f(q^{m\mu})=
{q^{-f(m+1)^2}+q^{-f(m-1)^2}-2q^{-fm^2}\over{2q^{-f}-2}}\>,$$ so
that $Z_K(M)=\phi_f(J'_{L_0}(\mu))$.

\vskip2ex\noindent{\bf Lemma 1}{\em \quad For any $a,b\in\ints$,
$l\in\nats$,
$\phi\big(q^{a\mu}(q^{\mu-b-1}-1)\ldots(q^{\mu-b-l}-1)\big)$ is
divisible by the $l\th$ cyclotomic polynomial $\phi_l(q)$ in
$\ints[q^{\pm1}]$.}\vskip2ex

\noindent{\bf Proof} \quad Suppose $\zeta$ is a root of unity of
order $l$.  It is sufficient to show that the expression in the
lemma vanishes when evaluated at any such $q=\zeta$.  Expanding out
the product,
$$\phi\big(q^{a\mu}(q^{\mu-b-1}-1)\ldots(q^{\mu-b-l}-1)\big)
=\sum_{S\subset\{1,\ldots,l\}}(-1)^{l-|S|}q^{\Sigma{}S-b|S|}\phi(q^{(a+|S|)\mu})$$
where $\Sigma{}S$ and $|S|$ denote the sum of elements and the
number of elements in the set $S$, respectively. Under the map
$S\lmt{}S+1$, which adds 1 (modulo $l$) to every element of $S$, the
number $\Sigma{}S$ increments by $|S|$ while $|S|$ stays unchanged.
Under repeated application of this map, we get orbits of size $l$
(for $0<|S|<l$) and two orbits with one element, namely the empty
set and the whole set.  The contribution from any orbit of size $l$
contains a factor $\sum_{i=0}^{l-1}q^{i|S|}$ which vanishes when
evaluated at $q=\zeta$.  The contribution from the two size one
orbits is
$$(-1)^l\phi(q^{a\mu})+q^{{1\over2}l(l+1)-bl}\phi(q^{(a+l)\mu})\>.$$
When evaluated at $q=\zeta$, $\phi(q^{a\mu})=\phi(q^{(a+l)\mu})$
while $\zeta^l=1$ and $\zeta^{{1\over2}l(l+1)}=(-1)^{l+1}$, and
again the contribution vanishes.\hfill$\Box$

Now for any $p\leq{}l$,
$q^{a\mu}(q^{\mu-b-1}-1)\ldots(q^{\mu-b-l}-1)$ can be written as a
linear combination (over $\ints[q^{\pm1}]$) of expressions of a
similar form with chains of length $p$,
$q^{a'\mu}(q^{\mu-b'-1}-1)\ldots(q^{\mu-b'-p}-1)$.  As a corollary,
the expression in Lemma~1 is divisible by $\phi_p(q)$ for any
$p\leq{}l$ and therefore also (since they are coprime) by
$<l>!\equiv\prod_{p=1}^l\phi_p(q)$. Here $<l>!\in\ints[q]$ is the
minimal monic polynomial in $q$ which vanishes at all roots of unity
of order $\leq{}l$ (but not at $q=1$).

\vskip2ex\noindent{\bf Corollary}{\em \quad For any $a,b\in\ints$,
$l\in\nats$,
$\phi\big(q^{a\mu}(q^{\mu-b-1}-1)\ldots(q^{\mu-b-l}-1)\big)$ is
divisible by $<l>!$ in $\ints[q^{\pm1}]$.}\vskip2ex

\vskip2ex\noindent{\bf Lemma 2}{\em\quad For any integers $a$,
$b_j$,
 $\phi_f\bigg(q^{a\mu}\prod_{j=1}^l(q^{\mu-b_j}-1)\bigg)$
 is divisible by $h^{\lceil{l\over2}\rceil}$.} \vskip2ex

\noindent{\bf Proof} \quad Let $\psi\takes{}B\lra\ints[[h]]$ be
defined by $\psi(q^{m\mu})=q^{m^2}$ for $m\in\ints$, extended
linearly over $\ints[[h]]$.  Write $q=e^\hbar$ and extend $\psi$
linearly to formal power series with rational coefficients. Then
$\psi(e^{m\mu\hbar})=e^{m^2\hbar}$. Comparing coefficients of powers
of $m$, $\mu^{2n}\lmt{(2n)!\over{}n!\hbar^n}$ while
$\mu^{2n+1}\lmt0$. In particular, for $s\leq{}r$,
$\psi(\hbar^r\mu^s)$ is divisible by $h^{\lceil{r/2}\rceil}$.

Now $\phi_{-1}$ and $\psi$ are related by
$$\phi(f)={1\over2h}\psi\big((q^\mu+q^{-\mu}-2)f\big)
={1\over2h}\psi\big((q^\mu-1)(q^{-\mu}-1)f\big)$$
 Expanding as a power series in $h$, whose coefficients are
polynomials in $\mu$,
$q^{a\mu}\prod_{j=1}^l(q^{\mu-b_j}-1)(q^\mu-1)(q^{-\mu}-1)$ is
divisible by $h^{l+2}$, and such that the coefficient of $h^r$
($r\geq{}l+2$) is a polynomial in $\mu$ of degree at most $r$.
Replacing $h$ by $\hbar$ preserves these properties (though
coefficients of powers of $\hbar$ may now be polynomials with
non-integral, but rational, coefficients). Hence applying $\psi$ we
obtain an element of $\rats[[\hbar]]$ divisible by
$\hbar^{\lceil{l+2\over2}\rceil}$. The lemma now follows.
\hfill$\Box$

\vskip2ex\noindent{\bf Lemma 3}{\em\quad For any integers $a$,
$b_j$, $l_j$,
 $$\phi_f\bigg(q^{a\mu}\prod_j(q^{\mu-b_j-1}-1)\cdots
 (q^{\mu-b_j-l_j}-1)\bigg)$$
 is divisible by $h^{\lceil{1\over2}\sum_jl_j\rceil}<\max{l_j}>!$
in $\ints[[h]]$ and in $\ints[q^{\pm1}]$.} \vskip2ex

\noindent{\bf Proof} \quad Note that $\phi_1$ and $\phi_{-1}$
interchange when $q$ and $q^{-1}$ are interchanged, so it is
sufficient to prove the lemma for $\phi\equiv\phi_{-1}$. Next
observe that the expression in the lemma is actually a polynomial in
$q^{\pm1}$ and so divisibility in $\ints[[h]]$ and in
$\ints[q^{\pm1}]$ are equivalent. Since $h$ and $<n>!$ are coprime
in $\ints[[h]]$, it is sufficient to prove divisibility by
$h^{\lceil{1\over2}\sum_jl_j\rceil}$ and by $<\max{l_j}>!$
separately; for the first we work in $\ints[[h]]$ (Lemma~2) and the
second in $\ints[q^{\pm1}]$ (Corollary).\hfill$\Box$

Putting Lemma~3 together with Theorem~1 from the last section gives
the following.

 \noindent{\bf
Theorem~2}{\em \quad The WRT invariant of a manifold, $M$, obtained
by $f=\pm1$ surgery around a knot, $L$,  is given by
$$Z_K(M)=\kern-10pt\sum_{l_j=0\atop{}j=1,\ldots,c}^\infty
 q^{c^0(\bfl)}
 \phi_f\bigg(q^{d(\bfl)\mu}
   \prod_{j=1}^c\prod_{p=1}^{l_j}(q^{\sigma_j(\mu-s(j)-p)}-1)\bigg)
   \kern-4pt\prod_{j=1}^c \binomials{s(j)+l_j}{l_j}_{\sigma_j}
$$
in which the contribution from $\bfl$ is divisible by
$h^{\lceil{1\over2}\sum_jl_j\rceil}<\max{l_j}>!$.  Here
$s(j)=-\sum_{i<j-}\epsilon_il_i$ while $c^0(\bfl)$ and $d(\bfl)$ are
given by (\ref{cl}) and (\ref{dl}) respectively. This is meaningful
when $q$ is a root of unity (then the sum is finite as the only
non-vanishing contributions are when $l_j<K$ for all $j$; the sum is
the WRT invariant) and as a formal power series in $h$ (the sum is
the Ohtsuki series, $Z_\infty(M)$).}\vskip2ex

In \cite{habiro}, Habiro obtains a more symmetric cyclotomic form
for Ohtsuki series which is needed to prove integrality in general.
For the case of surgery around knots, the above simple-minded
approach is sufficient.

\subsection{Growth rates in Ohtsuki series}

Suppose that $M$ is an integral homology sphere obtained by surgery
around a knot.  Then the Ohtsuki series $Z_\infty(M)$ is a formal
power series $\sum_{n=0}^\infty\lambda_n(M)h^n$ in $h=q-1$. In this
section we use the formula given in Theorem~3 to determine a bound
on the growth rate of coefficients $\lambda_n(M)$ with $n$.  By the
theorem, the only terms $\bfl$ which contribute to $\lambda_N(M)$
are those for which $\lceil{1\over2}\sum_jl_j\rceil\leq{}N$, that is
$\sum_jl_j\leq2N$; there are $\bin{c+2N}{c}=O(N^c)$ such $\bfl$'s.

The double product in the summand in Theorem~2 is a product of
$\sum_jl_j$ factors of the form $q^{\pm\mu+a}-1$ where $a\in\ints$.
Thus the contribution from $\bfl$ to $Z_\infty(M)$ can be expanded
out to give a sum of $2^{\sum_jl_j}\leq4^N$ terms of the form
 \begin{eqnarray}\label{zinfterm}
  \pm\phi_f(q^{\alpha\mu+\beta})
  .\prod_j\binomials{s(j)+l_j}{l_j}_{\sigma_j}
 \end{eqnarray}
 where
$\alpha=d(\bfl)+\sum_j\sigma_j|P_j|$ and
$\beta=c^0(\bfl)-\sum_j\sigma_j\big(s(j)|P_j|+\Sigma(P_j)\big)$,
while $P_j\subset\{1,\ldots,l_j\}$. Using (\ref{dl}) and (\ref{cl}),
one finds estimates
$$|\alpha|\leq2N(1+\max{|q_j|})+{b\over2},
|\beta|\leq(16c+10)N^2+2N(2+\max|q_j|)+{b\over2}$$ that is,
$\alpha=O(cN)$ and $\beta=O(cN^2)$, so that the order of the
coefficient of $h^n$ ($n\leq{}N$) in $\phi_f(q^{\alpha\mu+\beta})$
is that of $\bin{ac^2N^2}{n+1}$ ($a$ is a constant).

The final term to estimate is of the coefficient of $h^{N-n}$ in the
product of $c$ $q$-binomial coefficients on the right hand side of
(\ref{zinfterm}). Each of these $q$-binomial coefficients is of the
form $\binomials{a}{b}$ with $b\leq{}a\leq2N$; such a coefficient is
a sum of $\bin{a}{b}$ powers of $q$ (with exponents at most
$b(a-b)$) and so the coefficient of a power $h^r$ is at most
$$\bin{2N}{N}\bin{N^2}{r}\leq4^N\bin{N^2}{r}$$
Taking a product of $c$ such terms can accumulate in the coefficient
of $h^r$ at most $4^{cN}\bin{cN^2}{r}$ and thus the coefficient of
$h^N$ in (\ref{zinfterm}) is at most
 $$\sum_{n=0}^N\bin{ac^2N^2}{n+1}\bin{cN^2}{N-n}=\bin{(ac^2+c)N^2}{N+1}$$
Putting all the above estimates together gives the order of
$\lambda_N(M)$ to be at most that of
$$N^c.4^N\bin{ac^2N^2}{N+1}=x_N$$
which is an order of growth of the bound for which
$x_N\over{}x_{N-1}$ asymptotically grows linearly with $N$, with
slope of order $c^2$.

\noindent{\bf Theorem~3}{\em \quad There exist bounds $x_n$
dependent only on $n$ and the number of crossings $c$ in the knot
$K$, for which $|\lambda_N(S^3_K)|<x_N$ for all $N\in\nats$ while
$\lim\limits_{n\rightarrow\infty}{x_{N+1}\over{}Nx_N}$ exists and is
of order $c^2$.}\vskip2ex

This should be compared with the slope conjecture of
\cite{jacobylawrence} which states that
$\lambda_N(M)\over\lambda_{N-1}(M)$ is asymptotically linear in $N$,
with coefficient, $\sigma(M)$, called the slope of $M$. Combining
with the above result, it is reasonable to conjecture that
$\sigma(M)$ is of order $c^2$, when $M$ is obtained by surgery
around a knot with $c$ crossings. This is consistent with known
results (\cite{jacobylawrence}, \cite{lawrencerozansky}) on the
slope for Seifert fibred manifolds.

\section{Examples}

In the literature (see \cite{garoufalidisle}, \cite{habirocalc},
\cite{hikami}, \cite{letwobridge}, \cite{masbaum}) there are
computations of the colored Jones function for the simplest knots,
the trefoil, figure-8 knot and more generally, torus knots and twist
knots. Applying Theorems~1 and 2, we reproduce the simplest of these
results, although some are in different forms.

\subsection{The trefoil}

Suppose that $L$ is the right hand trefoil, presented as in
Figure~6. In general, the formulae obtained will be simpler if the
basepoint is chosen so that the first crossing encountered is an
overcrossing (and the last before returning to the basepoint is an
undercrossing. Crossing 3 is a left pointing crossing and the
appropriate enhanced Gauss diagram is shown on  the right in
Figure~6.

\begin{figure}
  \includegraphics[width=400pt]{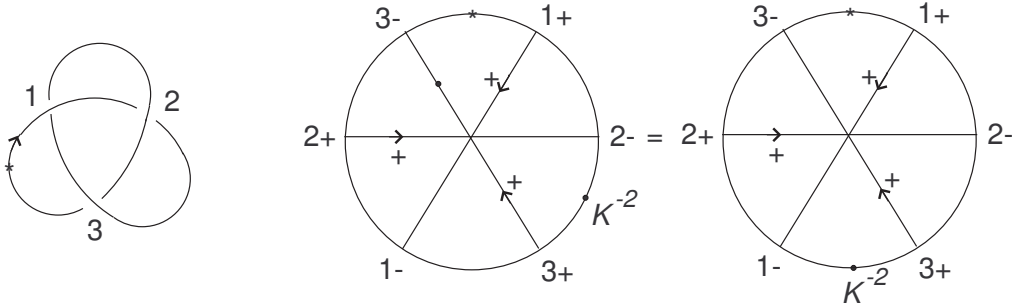}\\
  \caption{Trefoil knot and its enhanced Gauss diagram}
  \label{trefdiag}
\end{figure}

In this case $c=3$ and the six points (from crossings) on the outer
circle of the Gauss diagram are marked $1+$, $2-$, $3+$, $1-$, $2+$
and $3-$, respectively. Also there is just one blob $K^{-2}$, so
that $b=1$ and $\delta_1=-1$. The signs of the crossings $\sigma_j$
are all $+1$. Hence
$$s(1)=l_2-l_1-l_3\>,\quad
 s(2)=-l_1\>,\quad
 s(3)=-l_3\>.$$
 Non-vanishing contributions to $V'_{L_0}(\mu)$ appear only when
$s(j)$ are all non-negative, so that $l_1=l_3=0$. Then $s(1)=-l_2$
while $s(2)=s(3)=0$. We calculate $q_2=-2$ so that
 $$d(l)=-1-2l\>,\quad c(l)=1+{1\over2}l+{1\over2}l^2\>,$$
and so $c^\mu(l)=-(1+2l)\mu+(1+{1\over2}l+{1\over2}l^2)$. Theorem~1
now gives the Jones function of the trefoil as
$$J'_{L_0}(\mu)=\sum_{l=0}^\infty{}q^{-(2l+1)\mu+{l(l+1)\over2}+1}
 (q^{\mu-1}-1)\ldots(q^{\mu-l}-1)\>,$$
 while the Kashaev invariant is
 $$K_L(q)=\sum_{l=0}^\infty{}q^{{l(l+1)\over2}+1}
 (q^{-1}-1)\ldots(q^{-l}-1)\>.$$
 These formulae may be put in more standard form as
$$J'_{L_0}(\mu)=q\sum_{l=0}^\infty{}q^{-(l+1)\mu}(q^{1-\mu};q)_l\>,
\quad
 K_L(q)=q\sum_{l=0}^\infty(1-q)\ldots(1-q^l)\>,$$
 where $(a;q)_k=(1-a)(1-aq)\ldots(1-aq^{k-1})$.
 This formula was obtained by Habiro\cite{habiro} and Le, as well as
 another one,
 $$J'_{L_0}(\mu)=\sum_{k-0}^\infty(-1)^kq^{-k(k+3)\over2}q^{-k\mu}
 (q^{\mu+1};q)_k(q^{\mu-1};q^{-1})_k$$

\subsection{The figure 8 knot}

The figure 8 knot ($4_1$) is naturally presented in blackboard
framing as a zero framed knot. Figure~5 shows its enhanced Gauss
diagram along with the labeling of points on the circle.  Here $c=4$
with signs $\sigma_1=\sigma_3=1$ and $\sigma_2=\sigma_4=-1$. There
are two blobs $K^{-2}$ so that $b=2$, $\delta_a=-1$. Similar to the
case of the trefoil, the only non-vanishing contributions come from
$l_1=l_3=0$. Set $l_2=l$ and $l_4=m$. Then $s(1)=s(4)=l-m$ and
$s(2)=s(3)=0$. Also $q_1=q_3=-1$, $q_2=q_4=1$ so that from
(\ref{dl})
$$d(l,m)=l+m+1\>.$$

The last term in (\ref{cl}) requires an enumeration of all triples
$(i,k,j)$ for which $i<j-$ and $k<j+$; the part of the list for
which $i,k\in\{2\pm,4\pm\}$ is $(2+,2-,3)$, $(2+,4+,3)$,
$(2-,2-,3)$, $(2-,2-,4)$, $(2-,4+,3)$, $(4+,2-,3)$, $(4+,2-,4)$,
$(4+,4+,3)$, $(4-,2-,3)$ and $(4-,4+,3)$. Substituting into
(\ref{cl}) gives
$$c^0(l,m)=-(l+m+1)-l+(-l+m)+{1\over2}(l+m)
 +{1\over2}(-l^2-m^2)+(-l^2+lm)$$
Theorem~1 supplies the Jones function of the knot $4_1$ as
$$J'_{L}(\mu)=\kern-4pt\sum_{l\geq{}m\geq0}^\infty{}q^{c^\mu(l,m)}
 (q^{1-\mu}-1)\ldots(q^{l-\mu}-1)
 \binomials{l}{m}_-(q^{l-\mu-m+1}-1)\ldots(q^{l-\mu}-1)$$
where
$c^\mu(l,m)=(l+m+1)\mu-{3\over2}l^2+lm-{1\over2}m^2-{5\over2}l+{1\over2}m-1$.
The Kashaev invariant is then
$$K_L(q)=\sum_{l\geq{}m\geq0}q^{c^0(l,m)}(q-1)\ldots(q^l-1)(q-1)\ldots(q^m-1)
\binomials{l}{m}_+\binomials{l}{m}_-$$ This is to be compared with
Habiro \cite{habiro} and Le's formulae
$$J'_L(\mu)=\sum_{k=0}^\infty{}q^{k\mu}(q^{-\mu-1};q^{-1})_k(q^{-\mu+1};q)_k\>,
\quad K_L(q)=\sum_{k=0}^\infty(q^{-1};q^{-1})_k(q;q)_k$$
 in which
the symmetry $q\leftrightarrow{}q^{-1}$ due to the amphichirality of
$4_1$ is apparent.

\vskip2ex \noindent{\bf Acknowledgements} {The first author would
like to thank Thang Le and Lev Rozansky  for many conversations on
the subject of computations of 3-manifold invariants.  The second
author would like to acknowledge partial support from an Excellence
Grant from the Hebrew University.}

\end{document}